\input amstex
\documentstyle{amsppt}

\hsize=4.75in
\vsize=8in
\rightheadtext {Cohomology and Injectivity}
\leftheadtext {Erik Christensen and Allan M. Sinclair}
\NoBlackBoxes
\topmatter 
\title 
A cohomological characterization of approximately\\ 
finite dimensional von Neumann algebras 
\endtitle 
\author 
Erik Christensen and Allan M. Sinclair 
\endauthor 
\affil
University of Copenhagen and University of Edinburgh
\endaffil
\endtopmatter 
\document 
\heading
1. Introduction
\endheading 
One of the purposes in the computation of cohomology groups is to establish
invariants which may be helpful in the classification of the objects under
consideration. In the theory of continuous Hochschild cohomology for operator
algebras R. V. Kadison and J. R.  Ringrose proved \cite{10} that for any
hyperfinite von Neumann algebra $\Cal M$ and any dual normal $\Cal M$-bimodule
$\Cal S$, all the continuous cohomology groups vanish. Based on his
fundamental paper \cite{4} on injective von Neumann algebras A. Connes
established in 1978 in \cite{5} a converse statement to the result by Kadison
and Ringrose. Thus we have a characterization of approximately finite
dimensionality for von Neumann algebras by the vanishing of all continuous
Hochschild cohomology. When examining the proofs in \cite{1, 5} one finds that
there is a single module and a single derivation which has the characterizing
property such that $\Cal M$ is an injective von Neumann algebra if and only if
this test derivation is inner. The equivalence between injectivity and
approximately finite dimensionality was first established in \cite{4}. A later
proof can be found in \cite{7} and \cite{6} contains an extension to algebras
on non separable Hilbert spaces.  Having this we can state that a von Neumann
algebra is approximately finite dimensional if and only if Connes' test
derivation or test 1-cocycle is a coboundary. Various researchers have asked
whether the vanishing of every continuous second cohomology group for a von
Neumann algebra $\Cal M$ with coefficients in a dual normal module would imply
that $\Cal M$ is approximately finite dimensional. We can not settle this
question here, but we are able to construct a Banach $\Cal M$-bimodule $\Cal
S$ and a continuous 2-cocycle $\Phi$ on $\Cal M$ with coefficients in $\Cal S$
such that $\Cal M$ is approximately finite dimensional if and only if $\Phi$
cobounds. This is the only result in this paper. The module $\Cal S$ and the
2-cocycle $\Phi$ is constructed as a ``2-dimensional'' straightforward
generalization of the module and the test derivation from \cite{5}. In
\cite{5} the derivation $\delta:\Cal M\ @>>> \Cal Y\subseteq B(B(H))$ is given
by $\delta(m)(x)=[m,x]$. The exact description of the $\Cal M$-module $\Cal Y$
is left out here, but $\Cal Y$ is clearly a dual normal $\Cal M$-bimodule.
The 2-cocycle $\Phi$ we use in this example is constructed analogously by
 $$ 
\forall m_1,m_2\in\Cal M\;\;\forall x,y\in B(H): 
\Phi(m_1,m_2)(x,y)=[m_1,x][m_2,y]\;. 
 $$ 
Clearly $\Phi$ becomes a bilinear operator from $\Cal M\times\Cal M$ to 
the space of bilinear operators on $B(H)\times B(H)$ with values in 
$B(H)$. The module we choose is quite reasonable and the test 
2-cocycle is the obvious one. But we must for technical reasons make 
some assumptions on the module which prevents the module from being 
dual. The problem we encounter is that we want to ``lift'' some 
derivations simultaneously in a linear and continuous fashion. The 
only tool we could find to use is the combined results by Johnson \& 
Parrott \cite{9} and Popa \cite{12} on derivations of von Neumann 
algebras into the compact operators. In order to get ``sufficiently'' 
many derivations with values in the compacts the ``obvious'' module 
has to be reduced in such a way that it is no longer a dual Banach 
space. 
 
\medskip 
 
As general reference books we use \cite{11} for operator algebra theory and
\cite{8, 13} for  the theory of bounded Hochschild cohomology.
 
\heading
2. The result
\endheading 
 
We will now describe the test cocycle in detail and give the proof. 
Hence we consider in the following a von Neumann algebra $\Cal M$ on 
a Hilbert space $H$ and we define the $\Cal M$-bimodule $\Cal S$ 
as a subspace of the completely bounded bilinear operators on $B(H)$. 
Recall from \cite{2} that a bounded bilinear operator $F:B(H)\times 
B(H) @>>> B(H)$ can be extended to a bounded bilinear operator 
$F_n:M_n(B(H))\times M_n(B(H)) @>>> M_n(B(H))$ by the convention 
$F_n((a_{ij}),(b_{ij}))_{st}=\sum\limits_{i}F(a_{si},b_{it})$. If 
$\sup\limits_{n}\Vert F_n\Vert<\infty$, then we say that $F$ is 
completely bounded and define the completely bounded norm of $F$, 
$\Vert F\Vert_{cb}=\sup\limits_{n}\Vert F_n\Vert$. The space of 
completely bounded bilinear operators on $B(H)$ is denoted 
$CB(B(H)\times B(H),\;B(H))$. This space is a two-sided Banach $\Cal 
M$ bimodule by the elementary products, $(m. F)(x,y)=mF(x,y)$ $(F . 
m)(x,y)=F(x,y)m$. From the following description of the subspace 
$\Cal S$ it is easily seen that $\Cal S$ is also a Banach $\Cal M$ 
bimodule under this product. In the following $K(H)$ denotes the 
compact operators on $H$.   
  
\subheading{2.1 Description of the module}

The module $\Cal S$ is given as the linear space 
$$ 
\aligned 
\Cal S =  \lbrace & F\in CB(B(H)\times B(H),B(H))\mid F|\Cal M'\times 
B(H)=0\;\;\text{and}\\ 
& F|B(H)\times\Bbb C I=0\;\;\text{and}\;\;F:B(H)\times K(H) @>>> 
K(H)\;\;\text{and}\\ 
& F\;\;\text{is separately ultraweakly continuous in the right}\\ 
& \text{(second) variable}\rbrace\;, 
\endaligned 
 $$ 
equipped with the completely bounded norm $\Vert\quad\Vert_{cb}$.

\smallskip
 
The last two conditions are included in order to be able to use the 
existing literature on derivations \cite{9, 12} associated with the 
space of compact operators as a module. On the other hand 
these conditions make $\Cal S$ a non dual module. The reason why dual 
modules are preferred is that pathologies may occur in the non dual 
case. It is our judgment that $\Cal S$ is a reasonably large and 
natural module, such that the following example is healthy rather 
than pathological.

\subheading{2.2 Description of the 2-cocycle}

The mapping $\Phi:\Cal M\times\Cal M @>>> \Cal S$ is given by 
 $$ 
\forall m_1,m_2\in\Cal M\;\;\forall x,y\in B(H): 
\Phi(m_1,m_2)(x,y)=[m_1,x][m_2,y]\;. 
 $$ 
 It is immediate to see that $\Phi$ is a bilinear operator on $\Cal 
M$ with values in $\Cal S$. By \cite{8, 13} the coboundary of $\Phi$ 
is given by 
 $$ 
\aligned 
& \Delta\Phi(m_1,m_2,m_3)(x,y)=m_1[m_2,x][m_3,y]-\\ 
& [m_1m_2,x][m_3,y]+[m_1,x][m_2m_3,y]-[m_1,x][m_2,y]m_3= 0 
\endaligned 
 $$ 
 and we get that $\Phi$ is a 2-cocycle, i\.e\. $\Phi\in\Bbb 
Z_c^2(\Cal M,\Cal S)$. 
 
We will say that $\Phi$ is a coboundary if there exists a continuous 
linear mapping $\psi:\Cal M @>>> \Cal S$ such that $\Phi=\Delta\psi$, 
i\.e\. 
 $$ 
\aligned 
& \forall m_1,m_2\in\Cal M\;\;\forall x,y\in B(H):\\ 
& \Phi(m_1,m_2)(x,y)=(m_1\psi(m_2)-\psi(m_1m_2)+\psi(m_1)m_2)(x,y)\;. 
\endaligned 
 $$

 \proclaim{2.3 Theorem}
 The von Neumann algebra $\Cal M$ is approximately finite-dimen- sional if and
only if $\Phi$ is a coboundary.   
 \endproclaim 
 
\demo{Proof} In both directions of the proof we need the following 
observation regarding $\Phi$. Let $\rho:\Cal M @>>> CB(B(H)\times 
B(H),B(H))$ be defined by 
 $$ 
\forall m\in\Cal M\;\;\forall x,y\in B(H):\rho(m)(x,y)=[x,m]y\;. 
 $$ 
 Then the following computation shows that with respect to the $\Cal 
M$-bimodule $CB(B(H)\times B(H),B(H))$ we have $\Delta\rho=\Phi$. In fact 
for $m_1,m_2$ in $\Cal M$ and $x,y$ in $B(H)$ 
 $$ 
\aligned 
\Delta\rho(m_1,m_2)(x,y) & = 
m_1[x,m_2]y-[x,m_1m_2]y+[x,m_1]ym_2\\ 
& = [x,m_1](-m_2y+ym_2)=\Phi(m_1,m_2)(x,y)\;. 
\endaligned 
 $$ 
 By a similar argument we get that for $\gamma:\Cal M @>>> 
CB(B(H)\times B(H),B(H))$ given by $\gamma(m)(x,y)=x[m,y]$ 
we have $\Delta\gamma=\Phi$ 
 
Let us first suppose that $\Cal M$ is approximately 
finite-dimensional, $AF$, then $\Cal M'$ is  too and there exists an 
$\Cal M'$-bimodule conditional expectation $\pi$ of $B(H)$ onto $\Cal 
M'$. Let $\widehat{\gamma}:\Cal M @>>> CB(B(H)\times B(H),B(H))$ be 
given by 
 $$ 
\widehat{\gamma}(m)(x,y)=(x-\pi(x))[m,y] 
 $$ 
 then $\widehat{\gamma}(m)$ is in $\Cal S$ and by the previous 
computations 
 $$ 
\aligned 
\Delta\widehat{\gamma}(m_1,m_2) & = 
\Delta\gamma(m_1,m_2)(x-\pi(x),y)\\ 
& = \Phi(m_1,m_2)(x,y)\;. 
\endaligned 
 $$ 
 To prove the other implication we show in the following lines that 
if $\Phi$ cobounds then there must necessarily be a conditional 
expectation onto $\Cal M'$ and hence $\Cal M'$ and also $\Cal M$ must be 
injective and then approximately finite-dimensional. Suppose 
$\varphi:\Cal M @>>> \Cal S$ is a continuous linear mapping such that 
$\Delta\varphi=\Phi$. 
 
Let $P$ in $\Cal M$ denote the central projection corresponding to 
the continuous central summand of $\Cal M$ (i\.e\. $(I-P)\Cal M$ is 
of type I and $P\Cal M$ has no central summand of type I). As it is 
well known $(I-P)\Cal M$, (the type I part), is injective so we will 
focus on the algebra $P\Cal M$ on $PH$. It is easy to see that when 
$\Phi$ is restricted to this setting, $\Cal S$ is modified 
correspondingly and $\widehat{\varphi}$ is defined by 
 $$ 
\forall m_1\in P\Cal M\;\;\forall x,y\in B(PH)\;\; 
\widehat{\varphi}(m_1)(x,y):=P\varphi(m_1P)(xP,yP)P 
 $$ 
 then the coboundary of $\widehat{\varphi}$ equals the restricted 
$\Phi$. Hence we may and will in the following arguments assume that 
$\Cal M$ is a continuous von Neumann algebra. The basic fact we want 
to use in the following arguments is the combined results of \cite{9} 
and \cite{12} which can be stated as follows: 
 
\medskip 
 
Let $\Cal M$ be a continuous von Neumann algebra on a Hilbert space 
$H$ and let $\delta:\Cal M @>>> K(H)$ be a derivation. Then there is 
exactly one compact operator $k_\delta$ such that 
 $$ 
\forall m\in\Cal M:\delta(m)=[k_\delta,m]\;,\quad\text{and}\quad 
\Vert k_\delta\Vert\le\Vert\delta\Vert\;. 
 $$ 
 The norm estimate is not contained explicitly in the papers \cite{9, 
12}, but it can be deduced as follows. Let $\delta$ be a derivation 
of $\Cal M$ into the compacts $K(H)$, then by \cite{12} we know that 
there exists a unique $k$ in $K(H)$ such that $\delta(m)=[k,m]$. The 
uniqueness comes from the assumption that $\Cal M$ is a continuous 
von Neumann algebra. On the other hand this assumption implies that 
$\Cal M$ has a maximal Abelian von Neumann subalgebra $\Cal A$ which 
is diffuse, i\.e\. without minimal projections. Let $\Cal C$ denote 
the ultraweakly closed convex hull of the set $\{u\delta(u^*)\mid u 
\;\;\text{unitary in}\;\;\Cal A\}$, then by the proof of \cite{9, 
Th. 2.1} we find that there exists a $c\in\Cal C$ such that 
$\delta(a)=[-c,a]$ for $a$ in $\Cal A$ and moreover $c\in K(H)$. By 
assumption  
$$ 
\forall a\in\Cal A:[k,a]=\delta(a)=[-c,a] 
$$ 
so $(k+c)$ belongs to $\Cal A'\cap K(H)=0$ since $\Cal A$ is diffuse. 
By construction we have $\Vert c\Vert\le\Vert\delta|\Cal 
A\Vert\le\Vert\delta\Vert$ so $\Vert k\Vert\le\Vert\delta\Vert$.  
  
\medskip 
 
By assumption $\Delta\varphi=\Delta\rho=\Phi$ so 
$\delta=\varphi-\rho$ must be a derivation of $\Cal M$ into 
$CB(B(H)\times B(H),B(H))$. In symbols we get 
 $\forall x,y\in\Cal B(H)\forall m_1,m_2\in\Cal M:$
$$ 0=m_1\delta(m_2)(x,y)- 
\delta(m_1m_2)(x,y)+\delta(m_1)(x,y)m_2 
 $$ 
 and for any fixed pair $x,y$ we then have a derivation $d_{x,y}:\Cal 
M @>>> B(H)$ given by 
 $$ 
\forall x,y\in\Cal B(H)\;\;\forall m\in\Cal M: 
d_{x,y}(m)=\delta(m)(x,y)\;. 
 $$ 
 Let now $(p_i)_{i\in J}$ be an upwards filtering universal net of finite rank 
self-adjoint projections in $B(H)$ such that $p_i\uparrow I$. Then by 
assumptions on $\varphi$ and definition of $\rho$ we have for $x,y\in 
B(H)$, $i\in J$ that $d_{x,p_i y}$ is a derivation of $\Cal M$ into 
$K(H)$ such that $\Vert d_{x,p_i y}\Vert\le(\Vert\varphi\Vert+2)\Vert 
x\Vert\;\Vert y\Vert$. Hence there exists a unique $k(x,p_i y)$ in 
$K(H)$ such that $\Vert k(x,p_i y)\Vert\le 
(\Vert\varphi\Vert+2)\Vert x\Vert\;\Vert y\Vert$ and for $m\in M$ 
 $$ 
(\varphi-\rho)(m)(x,p_i y)=[k(x,p_i y),m]\;.\tag1 
 $$ 
 By uniqueness it is clear that $k(x,p_i y)$ is bilinear in $x$ and 
$y$.  
 
Since $(p_i)_{i\in J}$ was chosen as a universal net and we have that 
$\{\Vert k(x,p_i y)\Vert\mid i\in J\}$ is a bounded set for fixed $x$ 
and $y$ the net $(k(x,p_i y))_{i\in J}$ must be weakly convergent say 
with limit $k(x,y)$. By definition $\varphi(m)(x,y)$ is separately 
ultraweakly continuous in the $y$-variable and $\rho(m)$ has the same 
property by construction so we get from (1) 
 $$ 
\align 
& (\varphi-\rho)(m)(x,y)=[k(x,y),m]\tag2\\ 
& k:B(H)\times B(H) @>>> B(H)\;\;\text{is bilinear}\tag3\\ 
& \text{and continuous}\;. 
\endalign 
 $$ 
For $y=I$ we get, since $\varphi$ has range in $\Cal S$   that $-\rho(m)(x,I)=[k(x,I),m]$, so
$[x,m]+[k(x,I),m]=0$  
and we can define a continuous mapping $\pi:B(H) @>>> \Cal M'$ by 
 $$  
\pi(x)=x+k(x,I)\;. 
 $$ 
 For $x$ in $\Cal M'$ we have $(\varphi-\rho)(m)(x,p_i)=0$ for all 
$i\in J$ so $k(x,I)=0$ and $\pi$ is a continuous projection onto 
$M'$. For any $n\in\Bbb N$ we can repeat the arguments given above 
with respect to $\varphi_n(m)$ and $\rho_n(m):M_n(B(H))\times 
M_n(B(H)) @>>> M_n(B(H))$. The uniqueness of the implementing compact 
operators ensure that $\Vert k_n\Vert$ is bounded by the inequality 
 $$ 
\forall x,y\in M_n(B(H)):\Vert k_n(x,y)\Vert\le 
(\Vert\varphi_n\Vert+2)\Vert x\Vert\;\Vert y\Vert\;.\tag4 
 $$ 
 By definition of the norm on the module $\Cal S$   
we get that $\pi$ is completely bounded and 
 $$ 
\Vert\pi\Vert_{cb}\le 1+(\Vert\varphi\Vert+2)\;.\tag5 
 $$ 
 Following the result in \cite{3, Th\. 3.1} we then know that $\Cal 
M'$ is approximately finite-dimensional.  
 \enddemo 
\heading 
References
\endheading  
 
\item{[1]} J.~Bion-Nadal: {\it Banach bimodule associated to an 
action of a discrete group on a compact space}, Lecture Notes in 
Mathematics Vol.~{\bf 1132}, 30-37.
 
\item{[2]} E.~Christensen and A. M.~Sinclair: {\it A survey of 
completely bounded operators}, Bull. London Math. 
Soc. {\bf 21} (1989), 417-448.
 
\item{[3]} E.~Christensen and A. M.~Sinclair: {\it On von Neumann 
algebras which are complemented subspaces of $B(H)$}, J. 
Functional Analysis {\bf 122} (1994), 91-102.  
 
\item{[4]} A.~Connes: {\it Classification of injective 
factors}, Ann. Math. {\bf 104} (1976), 73-115. 
 
\item{[5]} A.~Connes: {\it On the cohomology of operator 
algebras}, J. Functional Analysis {\bf 28} (1978), 248-253. 
 
\item{[6]} G. A.~Elliott: {\it On approximately finite dimensional 
von Neumann algebras II}, Canad. Math. Bull. {\bf l21} (1978), 415-418. 
 
\item{[7]} U.~Haagerup: {\it A new proof of the equivalence of 
injectivity and hyperfiniteness for factors on a separable Hilbert 
space}, J. Functional Analysis {\bf l62} (1985), 160-201. 
 
\item{[8]} B. E.~Johnson: {\it Cohomology in Banach algebras},
Mem. Amer. Math. Soc. Vol.~{\bf 127} (1972). 
 
\item{[9]} B. E.~Johnson and S. K.~Parrott: {\it Operators 
commuting with a von Neumann algebra modulo the set of compact 
operators} J. Functional Analysis {\bf 11} (1972), 39-61. 
 
\item{[10]} R. V.~Kadison and J. R.~Ringrose: {\it Cohomology of 
operator algebras II. Extended cobounding and the hyperfinite 
case}, Arkiv f\"ur mathematik {\bf l9} (1971), 55-63. 
 
\item{[11]} R. V. Kadison and J. R. Ringrose: {\it Fundamentals of 
the theory of operator algebras}, Academic Press, New York 
1983.
 
\item{[12]} S. Popa: {\it The commutant modulo the set of compact 
operators of a von Neumann algebra}, J. Functional 
Analysis {\bf l71} (1987), 393-408. 
 
\item{[13]} A. M. Sinclair and R. R. Smith: {\it Hochschild cohomology of von
Neumann algebras}, London Mathematical Society Lecture Note Series, Vol.~{\bf
203}. 
\enddocument
\bye